\documentclass{article}
\usepackage[T2A]{fontenc}
\usepackage[cp1251]{inputenc}
\usepackage[english,russian]{babel}
\usepackage[tbtags]{amsmath}
\usepackage{amsfonts,amssymb,dsfont,babel,mathrsfs,amscd}

\newcommand{\abs}[1]{\left|#1\right|} 
\newcommand{\lmid}{\enskip\rule[-1ex]{0.4pt}{3ex}\enskip} 
\newtheorem{thrm}{Теорема}
\newtheorem{prpst}{Предложение}
\newtheorem{crlr}{Следствие}
\newtheorem{rmrq}{Замечание}

\newtheorem{bspl}{\hspace{15pt}Пример}
\begin{document}

{\LARGE
\begin{center}
О разложениях в ряды, связанных с дзета-функцией
\end{center}
}

{\Large\centerline{В.Е.Шестопал}}

\begin{abstract}

Using properties of the Riemann zeta-function we propose two new large classes of evaluated series.
Incidentally the first class represents integrals as generalized average on very nonuniform sequences. The second class contains inter alia a lot of new series with the Jacoby theta-functions and rationals of the exponential function. Moreover we propose many functions that can replace the Riemann zeta-function in similar constructions.

Two examples: 1) if $f(x)$ has period 1 and is in some Lipschitz class, we have for any natural $M>1$
$$
\ln M\cdot\int_0^1f(x)dx=\sum_{n\geq 1}
\sum_{k=1}^{M-1}\left[\frac{1}{Mn-k}f\left(\frac{\ln (Mn-k)}{\ln M}\right)-\frac{1}{Mn}
f\left(\frac{\ln (Mn)}{\ln M}\right)\right],
 $$

2) if $\varphi_{J,M,N}(w)=(-1)^J\left(\frac{d^J}{dw^J}\right)\left(\frac{N}{e^{Nw}-1}-\frac{M}{(e^{Mw}-1}\right),$ where $J,M,N$ are integer,\\ $M>N>1,$ $J\geq 0$ and for all
$n\in\mathds{Z}$
$\left(e^{M\left(M/N\right)^{n+w}}-1\right)\left(e^{N\left(M/N\right)^{n+w}}-1\right)\neq 0,$ we
have
$$
\sum_{n\in\mathds{Z}}(M/N)^{(J+1)(n+w)}\varphi_{J,M,N}((M/N)^{n+w})=J!.
 $$

На основе анализа дзета-функции Римана строятся два новые большие семейства точно
суммируемых рядов. Первое семейство используется для изучения свойств некоторых классов
неравномерно распределенных последовательностей, что позволяет получать интегралы от непрерывных функций как обобщенные средние значения по таким последовательностям. 
\end{abstract}

В работах Эйлера, Чебышёва, Римана, Адамара и де ля Валле-Пуссена раскрылась хорошо ныне
известная польза, которую приносят исследования свойств дзета-функции Римана в изучении распределений простых чисел \cite{V-K,T}. Большой интерес к единственной оставшейся недоказанной гипотезе Римана и ее многочисленные иные приложения в теории чисел послужили одним из мотивов предложенного Гильбертом поиска применений дзета-функции в других областях математики.

Теорема универсальности Воронина \cite{V-K} показывает особое место этой функции в комплексном анализе.

Главная цель данной работы --- полнее раскрыть возможности применения некоторых открытых Эйлером соотношений для дзета-функции в теории чисел и анализе, показывая способы построения с ее помощью обширных
семейств новых явно суммируемых рядов (первые примеры приведены в \cite{Sch}). 

Прежде всего это дает представления интегралов
$\int_0^1f(x)dx=\sum_{n>0}a_nf(b_n)$ в случае непрерывных, в том числе и весьма негладких, функций $f.$ В качестве последовательностей $b_n$ могут использоваться очень неравномерные последовательности, что отражается на коэффициентах $a_n.$ 
Эти ряды можно рассматривать как обобщенные средние значения по таким последовательностям.
Среди прочего это позволяет получить свойства последовательностей дробных долей логарифмов натуральных чисел, а также многих других связанных с их цепными дробями последовательностей. Другими применениями таких рядов служат разнообразные граничные представления голоморфных и гармонических функций.

Ряды иного типа показывают новые свойства экспоненты и тета-функций Якоби, причем эти две функции соответствуют лишь двум значениям вещественного параметра, характеризующего семейства новых точно суммируемых рядов.

В заключительной части работы указаны другие классы функций, свойства которых подобны свойствам дзета-функции, что позволяет получать сходными методами еще больше точно суммируемых рядов и соотношений между рядами.

\section{Обобщенные средние непрерывных функций по неравномерным последовательностям}
\hspace{13pt} Здесь будут приведены представления интегралов от периодических функций по
периоду в виде рядов. Первые возникающие при этих построениях разложения позволяют охарактеризовать свойства 
последовательностей вида $\{\log_Mn\lmid n\in\mathds{N}\}$ и их  дробных частей при (сначала целых, а затем) рациональных $M>1.$

В качестве одного из применений мы показываем, что в случае области с достаточно простой границей граничные
представления голоморфных функций области  могут быть получены с помощью рядов.

Введем несколько обозначений и определений:

1) $\mathds{N}$ --- множество натуральных (целых положительных) чисел, для всякого
вещественного $x$ величина
$\{x\}$ --- его дробная часть; дополнительно положим $\{x/0\}=0;$

2) $\mathbf{B}_{L},$ где $L\geq 0,$ --- множество положительных монотонных отображений $(0,+\infty)$ в
себя, для каждого $\rho$ из которых ряды $\sum_{n\geq 2}\rho(c/n^2)(\ln{n})^L/n$ сходятся
при $c>0;$

3) $\mathbf{P_B}_{L}$ --- множество функций на $\mathds{R}$ с периодом единица, каждой $h$
из которых можно сопоставить функцию $\rho\in\textbf{B}_L,$ удовлетворяющую требованию: $(f(x)-f(y))/\rho(\abs{x-y})$ непрерывна;

4) $\mathbf{C_B}_L$ --- пространство функций на единичной окружности, каждой $F$ из
которых можно сопоставить функцию $\rho\in\textbf{B}_L,$ удовлетворяющую требованию:
$(F(x)-F(y))/\rho(\abs{x-y})$ непрерывна;

5) $\mathbf{A_B}_L$ --- пространство функций аналитических внутри
$\mathbf{D}=\{z\lmid\abs{z}\leq 1\},$ и ограничения которых на единичную окружность
принадлежат $\mathbf{C_B}_L;$

6) будем сокращать $\mathbf{B}_0,\,\mathbf{P_B}_0,\,\mathbf{C_B}_0,\,\mathbf{A_B}_1$
соответственно до $\mathbf{B},\,\mathbf{P_B},\,\mathbf{C_B},\,\mathbf{A_B};$

7) $\Big\{x\Big\}$ --- дробная часть вещественного числа $x;$

8) если $x>0,$ то под $x^y$ всегда будет пониматься $e^{y\ln x}.$

\begin{prpst}
Пространство $\mathbf{P_B}$ включает в себя функции $h,$ удовлетворяющие условиям
Липшица, т.е $\abs{f(x)-f(y)}\leq C\abs{x-y}^a,$ $a>0,$ а также  функции с менее ограничительными требованиями
вида $\abs{f(x)-f(y)}\leq C\abs{\,\ln{\abs{\,\ln{\abs{x-y}}}}}^{-b}\abs{\,\ln{\abs{x-y}}}^{-1},$
где $b>1.$
\end{prpst}

\begin{prpst}
Если $f\in\mathbf{P_B}_L$ и $g\in\mathbf{P_B}_L,$ то $f+g\in\mathbf{P_B}_L\ni fg.$
\end{prpst}
\textbf{Доказательство}. Пусть $f(x)-f(y)=\rho_1(\abs{x-y})\varphi(x,y)$ и
$g(x)-g(y)=\sigma(\abs{x-y})\rho_2(x,y),$ где $\rho_1,\rho_2\in\mathbf{B}_L,$ а $\varphi,\chi$
ограничены. В качестве $\rho$ для $f+g$ и $fg$ можно взять $\rho_1+\rho_2.\Box$

\subsection{Первое представление}
\hspace{13pt} Начнем с применения одного из соотношений Эйлера для дзета-функции.
\begin{thrm}\label{Th-2}
Если целое $M>1,$ и $f(x)\in\mathbf{P_B},$ то 
\begin{equation}\label{eq-4}
\ln M\;\int_0^1f(x)dx=\sum_{n\geq 1}
\sum_{k=1}^{M-1}\left[\frac{1}{Mn-k}f\left(\frac{\ln (Mn-k)}{\ln M}\right)-\frac{1}{Mn}
f\left(\frac{\ln (Mn)}{\ln M}\right)\right].
\end{equation}

Ряд по $n$ сходится абсолютно.
\end{thrm}
\textbf{Доказательство.} Справедливо (особенно хорошо известное при $M=2$) соотношение
\begin{equation}\label{Hpt}
\left(1-M^{1-s}\right)\zeta(s)=\sum_{n\geq 1}\sum_{k=1}^{M-1}\left[
\frac{1}{(Mn-k)^s}-\frac{1}{(Mn)^s}\right],
\end{equation}
в правой части которого ряд равномерно сходится в области $\Re{(s)}>\varepsilon$ при любом 
$\varepsilon>0.$

Дзета-функция имеет в точке $s=1$ полюс с вычетом равным единице, а потому
$$
\sum_{n\geq 1}\sum_{k=1}^{M-1}\left[\frac{1}{Mn-k}-\frac{1}{Mn}\right]=\ln M.
 $$

Функция $1-M^{1-s}$ имеет нули вида $s_l=1+2\pi il/\ln M,$ $l\in\mathds{Z}.$ Следовательно
при целых $l$ справедливо равенство
$$
\sum_{n\geq 1}\sum_{k=1}^{M-1}\left[
\frac{1}{(Mn-k)^{1+2\pi il/\ln M}}-\frac{1}{(Mn)^{1+2\pi il/\ln M}}\right]=\delta_{l0}\ln M.
 $$
Тем самым в случае функций вида $f(x)=\sum_{l=-N}^Na_le^{2\pi i\,lx},$ соотношение
(\ref{eq-4}) доказано.

Переходя к более широкому классу $\mathbf{P_B},$ далее будем рассматривать правую часть
(\ref{eq-4}) как сумму $M-1$ рядов\\
\centerline{$\sum_{n\geq 1}\left[\frac{1}{Mn-k}f\left(\frac{\ln (Mn-k)}{\ln M}\right)-\frac{1}{Mn}
f\left(\frac{\ln (Mn)}{\ln M}\right)\right],$} исследуя для каждого из них сходимость и
аппроксимационные свойства.

Пусть
\begin{equation}\label{stm} 
\rho\in\mathbf{B}\quad\mbox{и}\quad f(x)-f(y)=\rho(\abs{x-y})\chi(x,y),
\end{equation}
где $\chi$ непрерывна и ограничена.

Прежде всего это делает очевидной сходимость ряда в правой части (\ref{eq-4}), поскольку
\begin{equation}\label{stmtn} 
\left|\frac{1}{Mn-k}f\left(\frac{\ln (Mn-k)}{\ln M}\right)-\frac{1}{Mn}f\left(\frac{\ln (Mn)}{\ln M}\right)\right|\leq
\end{equation}
$$
\leq\left|\frac{1}{Mn-k}\left( f\left(\frac{\ln (Mn-k)}{\ln M}\right)-f\left(\frac{\ln (Mn)}{\ln M}\right)\right)\right|+
\left|\left(\frac{1}{Mn-k}-\frac{1}{Mn}\right)f\left(\frac{\ln (Mn)}{\ln M}\right)\right|\leq
 $$
$$
\leq\frac{\sup{\abs{\chi}}}{Mn-k}\rho\Big(-\frac{\log_M(1-k/(Mn))}{n}\Big)+\frac{k\,\max{\abs{f}}}{Mn(Mn-k)}
 $$
и сходятся ряды $\sum_{n\geq 1}\rho(c/n^2)/n,$ $c>0.$

Построим специальную последовательность $f_J$ тригонометрических многочленов, которая
равномерно сходится к $f,$ а суммы рядов (\ref{eq-4}) для этих многочленов сходятся к сумме
для $f$ (необходимая сходимость интегралов в левой части (\ref{eq-4}) верна для любой равномерно сходящейся последовательности
тригонометрических многочленов). С этой целью можно воспользоваться ядром Фейера и
стандартной техникой использования свойств сингулярных ядер \cite{F}:
$$
f_J(x)=\frac{1}{J}\int_0^1f(x+t)\left(\frac{\sin{\pi Jt}}{\sin{\pi t}}\right)^2dt, \quad J\in\mathds{N}.
 $$ 

Необходимо доказать, что при $1\leq k\leq M-1$
$$
\sum_{n\geq 1}\left[\frac{1}{Mn-k}f_J\left(\frac{\ln (Mn-k)}{\ln M}\right)-\frac{1}{Mn}
f_J\left(\frac{\ln (Mn)}{\ln M}\right)\right] \rightarrow
 $$
$$
\rightarrow \sum_{n\geq 1}\left[\frac{1}{Mn-k}f\left(\frac{\ln (Mn-k)}{\ln M}\right)-\frac{1}{Mn}
f\left(\frac{\ln (Mn)}{\ln M}\right)\right].
 $$

Возьмем
$\delta\in(0,1/2)$ и запишем вариацию $f$ в виде 
$$
f(x)-f_J(x)=\frac{1}{J}\int_0^1(f(x)-f(x+t))\left(\frac{\sin{\pi Jt}}{\sin{\pi t}}\right)^2dt=
\frac{1}{J}\left(\int_{-\delta}^{+\delta}\ldots+\int_{+\delta}^{1-\delta}\ldots\right).
 $$

Введем сокращение:
$x_k=x_{k,n}=\log_M (Mn-k)$ и воспользуемся соотношением
$$
\frac{f(x_k)-f_J(x_k)}{Mn-k}-
\frac{f(x_0)-f_J(x_0)}{Mn}=
 $$
$$
=\left(\frac{1}{Mn-k}-\frac{1}{Mn}\right)\Big(f(x_k)-f_J(x_k)\Big)+\frac{1}{Mn}\Big(f(x_k)-f(x_0)-f_J(x_k) +f_J(x_0)\Big)=
 $$
\begin{equation}\label{ueq-1}
=\left(\frac{1}{Mn-k}-\frac{1}{Mn}\right)\Big(f(x_k)-f_J(x_k)\Big)+
\end{equation}
\begin{equation}\label{ueq-2}
+\frac{1}{Mn}\frac{1}{J}\left(\int_{-\delta}^{+\delta}\ldots+\int_{+\delta}^{1-\delta}\ldots\right)
\Big(f(x_k)-f(x_k+t)-f(x_0)+f(x_0+t)\Big)\left(\frac{\sin{\pi Jt}}{\sin{\pi t}}\right)^2dt.
\end{equation}

Выражение (\ref{ueq-1}) равномерно оценивается величиной $\abs{f(x)-f_J(x)}$ при $J\to +\infty,$
которая будет умножена на сумму сходящегося ряда. 
Абсолютная величина (\ref{ueq-2}) оценивается на основе неравенства\\
\centerline{$\abs{f(x_k)-f(x_0)-f(x_k+t)+f(x_0+t)}\leq \rho(\abs{x_k-x_0})\cdot
\abs{\chi(x_k,x_0)-\chi(x_k+t,x_0+t)},$}\\
что дает для (\ref{ueq-2}) верхнюю оценку
$$
\frac{\rho(\abs{x_{k,n}-x_{0,n}})}{Mn}\frac{1}{J}\left(\int_{-\delta}^{+\delta}\ldots+\int_{+\delta}^{1-\delta}\ldots
\right)\abs{\chi(x_k,x_0)-\chi(x_k+t,x_0+t)}\left(\frac{\sin{\pi Jt}}{\sin{\pi t}}\right)^2dt.
 $$
Ряд $\sum_{n\geq 1}\rho(\abs{x_{k,n}-x_{0,n}})/n=\sum_{n\geq 1}\rho\Big(-\frac{\log_M(1-k/(Mn))}{n}\Big)/n$ сходится. В отношении же интегралов  воспользуемся стандартными свойствами ядра Фейера \cite{F}: малость интеграла по
$[-\delta,\delta]$ определяется выбором $\delta$ и непрерывностью $\chi,$ а малость
интеграла по $[\delta,1-\delta]$ --- выбором достаточно большого $J.\Box$

\subsection{О последовательностях логарифмов и представлениях интегралов рядами}

\hspace{13pt} Доказанная Теорема \ref{Th-2}  дает любопытные свойства дробных долей
логарифмов натуральных чисел с рациональными основаниями.
Действительно, в \cite{K-M} показано, что последовательность дробных долей
$\{\log_An\},\,n\in\mathds{N},$ может быть весьма неравномерно распределенной на $[0,1]$
последовательностью, что рассматривается как одно из препятствий к ее применению в
численном интегрировании. Тем не менее из полученных представлений видно, что интегралы
от функций из $ \mathbf {P_B} $ можно заменять рядами по взвешенным значениям этих функций
на таких последовательностях. Соотношение (\ref{stmtn}) позволяет оценивать
сходимость ряда (\ref{eq-4}) в случае весьма негладкой функции $h.$

Имея Теорему \ref{Th-2} можно получить и другие представления интегралов в виде рядов по последовательностям "узлов суммирования"\, с подобными свойствами. Это же позволяет говорить об усреднении по многим другим неравномерным последовательностям.

\begin{thrm}
Пусть функция $f$ непрерывна, $\varphi$ непрерывно дифференцируема,
\\$\varphi(1)-\varphi(0)=1$ и $f(\varphi(x))\varphi'\Big(x\Big)\in\mathbf{P_B}.$

В этом случае
\begin{equation}\label{eq-04}
\ln M\int_0^1f(x)dx=\sum_{n\geq 1}
\sum_{k=1}^{M-1}\left[\frac{f\left(\varphi\Big(\log_M (Mn-k)\Big)\right)\varphi'\Big(\log_M (Mn-k)\Big)}{Mn-k}-\right. 
\end{equation}
$$
\left. -\frac{f\left(\varphi\Big(\log_M (Mn)\Big)\right)\varphi'\Big(\log_M (Mn)\Big)}{Mn}\right].
 $$
\end{thrm}
\textbf{Доказательство}. Имеем: $\int_0^1f(x)dx=\int_{\varphi(0)}^{\varphi(1)}f(\varphi(x))\varphi'(x)dx=\int_0^1f(\varphi(x))\varphi'(x)dx.\Box$

\begin{thrm}\label{Th-25} 
Пусть даны целые $M>1\leq L$ и дифференцируемая функция $\varphi,$ и пусть
$\chi(\varphi(x))=x,$ $\varphi(L)=1,$ $\varphi(0)=0.$

Если $f(x+1)=f(x)$ и
$f\Big(\chi\Big(\Big\{x\Big\}\Big)\Big)g\Big(\Big\{x\Big\}\Big)\in\mathbf{P_B}\,,$  то
$$
\ln M\!\int_0^1\!f\Big(\chi(x)\Big)g(x)dx\!=\!
\sum_{n\geq 1}\sum_{k=1}^{M-1}\left[\frac{f\!\left(\!\chi\Big(\Big\{\log_M(Mn-k)\!\Big\}\Big)\right)\!g\!\left(\Big\{\log_M(Mn-k)\Big\}\right)}{Mn-k}\,-\right. 
 $$
$$
\left. -\,\frac{f\left(\chi\Big(\Big\{\!\log_M(Mn)\!\Big\}\Big)\right)g\left(\Big\{\log_M(Mn)\Big\}\right)}{Mn}\right]=\ln\!M \int_0^1f(x)\sum_{l=0}^{L-1}g(\varphi(l+x))\varphi'(l+x)\,dx.
 $$
\end{thrm}
\textbf{Доказательство} основано на (1.) соотношениях
$$
\int_0^1\!f\Big(\chi(x)\Big)g(x)dx=\sum_{l=0}^{L-1}\int_{\varphi(l)}^{\varphi(l+1)}\ldots=\int_0^1f(x)\sum_{l=0}^{L-1}g(\varphi(l+x))\varphi'(l+x)\,dx
 $$
и (2.) последующей замене $f(x)$ на $f(x)/G_\varphi(x).\Box$

\begin{crlr}
Пусть $f(x+1)=f(x)$ и $G_\varphi(x)=\sum_{l=0}^{L-1}g(\varphi(l+\Big\{x\Big\}\,))\,\varphi'(l+\Big\{x\Big\}\,)$
причем\\ $f\Big(\chi\Big(\Big\{x\Big\}\Big)\Big)g\Big(\Big\{x\Big\}\Big)\Big/G_\varphi\Big(\Big\{\chi\Big(\Big\{x\Big\}\Big)\Big\}\Big)\in\mathbf{P_B}\,.$ В этом случае
$$
\ln M\, \int_0^1f(x)dx=\sum_{n\geq 1}\sum_{k=1}^{M-1}\left[\frac{f\!\left(\chi\Big(\Big\{\log_M(Mn-k)\!\Big\}\Big)\right)\!
g\!\left(\Big\{\log_M(Mn-k)\Big\}\right)}{(Mn-k)\,G_\varphi\left(\!\Big\{\chi\Big(\Big\{\log_M(Mn-k)\!\Big\}\Big)\Big\}\right)}\,-\right. 
 $$
\begin{equation}\label{eq-22} 
\left. -\,\frac{f\left(\chi\Big(\Big\{\!\log_M(Mn)\!\Big\}\Big)\right)g\left(\Big\{\log_M(Mn)\Big\}\right)}{Mn\,G_\varphi
\left(\Big\{\chi\Big(\Big\{\!\log_M(Mn)\!\Big\}\Big)\Big\}\right)}\right].
\end{equation}
\end{crlr}

\begin{bspl}
Положив $L=M$ и $\chi(x)=M^x-1,$ мы получим в качестве системы узлов суммирования
множество рациональных чисел вида $\Big\{p/M^q\Big\},$ где $p,q\in\mathds{N}.$
\end{bspl}

\begin{rmrq}
Используя в качестве исходного вместо разложения (\ref{eq-4}) Теоремы \ref{Th-2}
соотношения (\ref{eq-22}) (или подобные (\ref{eq-23})), можно неограниченное число раз
повторить все существенные моменты построений Теорем
\ref{Th-25} и \ref{Th-24}, при желании меняя на каждом шаге функции $g$ и получая
представления в виде рядов для интегралов от непрерывных функций из довольно широких
классов. Попутно мы будем получать всё новые семейства числовых последовательностей,
обеспечивающие эти представления.
\end{rmrq}

Следующее утверждение показывает возможность использования в качестве узлов
суммирования не только логарифмов натуральных чисел, но и величин, естественно
возникающих при их разложениях в непрерывные дроби. Одновременно оно является примером
обобщения предшествующего результата на случай $L=\infty.$
\begin{thrm}\label{Th-24}
Пусть даны целые $M>1\leq L.$

Если $f(x+1)=f(x)$ и $f\Big(L\Big/\Big\{x\Big\}\Big)g\Big(\Big\{x\Big\}\Big)\in\mathbf{P_B}\,,$
справедливы равенства
$$
\ln M\!\int_0^1\!f\Big(L\Big/x\Big)g(x)dx\!=\!
\sum_{n\geq 1}\sum_{k=1}^{M-1}\left[\frac{f\!\left(L\Big/\Big\{\log_M(Mn-k)\Big\}\right)\!
g\!\left(\Big\{\log_M(Mn-k)\Big\}\right)}{Mn-k}\,-\right. 
 $$
$$
\left. -\,\frac{f\left(L\Big/\Big\{\log_M(Mn)\Big\}\right)g\left(\Big\{\log_M(Mn)\Big\}\right)}{Mn}\right]=L\ln M\, \int_0^1f(x)\sum_{n\geq L}\frac{g(L/(n+x))}{(n+x)^2}\,dx.
 $$
\end{thrm}
\textbf{Доказательство} основано на соотношениях
$$
\int_0^1f\Big(L\Big/x\Big)g(x)dx=\sum_{n\geq L}\int_{L/(n+1)}^{L/n}f\Big(\Big\{L/x\Big\}\Big)g(x)dx=
L\int_0^1f(x)\sum_{n\geq L}\frac{g(L/(n+x))}{(n+x)^2}\,dx
 $$
и замене в п.2 $h$  на $h/G.\Box$

\begin{crlr}
Если функция $g(x)$ такова, что функция $G(x)$ может быть определена равномерно сходящимся
рядом $G(x)=\sum_{n\geq L}g\Big(L\Big/(n+\Big\{x\Big\}\Big)\Big)\Big(n+\Big\{x\Big\}\Big)^{-2},$
$f\Big(\Big\{x\Big\}\Big)=f(x)$ и\\
$f\Big(L\Big/\Big\{x\Big\}\Big)g\Big(\Big\{x\Big\}\Big)\Big/G\Big(L\Big/\Big\{x\Big\}\Big)\in
\mathbf{P_B},$ то верно равенство
$$
L\ln M\, \int_0^1f(x)dx\,=\,
\sum_{n\geq 1}\sum_{k=1}^{M-1}\left[\frac{f\left(L\Big/\Big\{\log_M(Mn-k)\Big\}\right)
g\left(\Big\{\log_M(Mn-k)\Big\}\right)}{(Mn-k)\cdot G\left(L\Big/\Big\{\log_M(Mn-k)\Big\}\right)}\,-\right. 
 $$
\begin{equation}\label{eq-23} 
\left. -\,\frac{f\left(L\Big/\Big\{\log_M(Mn)\Big\}\right)g\left(\Big\{\log_M(Mn)\Big\}\right)}{Mn\cdot G\left(L\Big/\Big\{\log_M(Mn)\Big\}\right)}\right].
\end{equation}
\end{crlr}

\begin{rmrq}
Итерируя процедуру получения представления интеграла рядом с помощью замены
$x\to\Big\{1\Big/\Big\{x\Big\}\Big\},$ мы получим формулы замены интегралов суммами
рядов, узлы суммирования которых --- полные частные \cite{B} логарифмов натуральных чисел, т.е. величины, возникающие на
заранее задаваемом шаге разложения чисел $\{log_Mn\in\mathds{N}\}$ в непрерывную дробь.
\end{rmrq}

\subsection{Связи с комплексным анализом}
\hspace{13pt} Предложенные выше ряды можно использовать использовать в
комплексном анализе. Например, применяя в Теореме \ref{Th-2} подстановку $f(x)=F(e^{2\pi ix}),$ получаем следующее утверждение.
\begin{thrm}\label{Th-2-an}
Если целое $M>1$ и $F(z)\in\mathbf{C_B},$ то 
\begin{equation}\label{eq-5-4}
\frac{\ln M}{2\pi i}\oint_{\abs{z}=1}\frac{F(z)}{z}dz= \sum_{n\geq 1}
\sum_{k=1}^{M-1}\Bigg[\frac{F\Big((Mn-k)^{2\pi i/\ln M}\Big)}{Mn-k}-
\frac{F\Big((Mn)^{2\pi i/\ln M}\Big)}{Mn}\Bigg].\Box
\end{equation}

Ряд по $(n,k)$ сходится абсолютно.
\end{thrm}

В частности большое семейство рядов из правой части (\ref{eq-5-4}) в случае $F(z)=G(z,\bar{z})$
с аналитической $G(v,w)$ суммируются при помощи теоремы Коши о вычетах \cite{M}. Эти же ряды можно использовать для получения с их помощью
представлений значений голоморфных в односвязных областях функций через их значения на границе. Вот
несколько примеров.

\begin{thrm}\label{Th-5}
1. Пусть целые $L\neq 0$  и $M>1,$ $\abs{a}=1$ и $f(z)\in\mathbf{A_B},$ в этом случае
\begin{equation}\label{eq-5}
\ln M\,f(0)=\sum_{n\geq 1}\sum_{k=1}^{M-1}\left[\frac{1}{Mn-k}f\Big(a(Mn-k)^{2\pi iL/\ln M}\Big)
-\frac{1}{Mn}f\Big(a(Mn)^{2\pi iL/\ln M}\Big)\right].
\end{equation}

2. Ели также $J\in\mathds{N}$ и $0<\abs{c}<1,$ то
$$
\sum_{n\geq 1}\sum_{k=1}^{M-1}\left[\frac{1}{Mn-k}\cdot\frac{f\Big((Mn-k)^{2\pi i/\ln M}\Big)}
{\abs{(Mn-k)^{2\pi i/\ln M}-c}^{2J}}
-\frac{1}{Mn}\cdot\frac{f\Big((Mn)^{2\pi i/\ln M}\Big)}{\abs{(Mn)^{2\pi i/\ln M}-c}^{2J}}\right]
 $$
$$
=\frac{\ln M}{(J-1)!}\,\frac{d^{J-1}}{dz^{J-1}}\,[z^{J-1}f(z)(1-\bar{c}z)^{-J}]\Big|_{z=c}\quad\mbox{и}
 $$ 
$$
\sum_{n\geq 1}\sum_{k=1}^{M-1}\left[\frac{1}{Mn-k}\cdot\frac{f\Big((Mn-k)^{2\pi i/\ln M}\Big)}
{\abs{(Mn-k)^{2\pi i/\ln M}-1/c}^{2J}}
-\frac{1}{Mn}\cdot\frac{f\Big((Mn)^{2\pi i/\ln M}\Big)}{\abs{(Mn)^{2\pi i/\ln M}-1/c}^{2J}}\right]=
 $$ 
$$
=\frac{\ln M}{(J-1)!}\,\abs{c}^{2L}\,\frac{d^{J-1}}{dz^{L-1}}\,[z^{J-1}f(z)(1-cz)^{-J}]
\Big|_{z=\bar{c}}\,.
 $$
\end{thrm}
\textbf{Доказательство}. 1. Соотношения (\ref{eq-5}) следуют из формулы (\ref{eq-5-4}) (и интегральной формулы
Коши), примененной к функциям $F(z)=f(az^L)$ (при $L>0$) или $F(z)=\overline{f(a\bar{z}^L)}$ (в случае $0$).

2. Следует из (\ref{eq-5-4}), примененной к
$$
F_1(z)=f(z)\abs{z-c}^{-2J}=f(z)(z-c)^{-J}(1/z-\bar{c})^{-J},
 $$
$$
F_2(z)=f(z)\abs{z-1/c}^{-2J}=f(z)(z-1/c)^{-J}(1/z-1/\bar{c})^{-J}.\Box
 $$

Соотношения (\ref{eq-5}) можно обобщить для получения граничных представлений аналитических функций в случае
более общих областей, с помощью очень разнообразных наборов граничных точек и соответствующих весов.

\begin{thrm}\label{Th-3}
Пусть целое $M>1,$  функция $\mu(z)f(g(z))$ принадлежит $\mathbf{A_B}$ и определена в нуле, причем $\mu(0)=1$.

Если  $g(0)=w,$ то
\begin{equation}\label{eq-6}
f(w)\ln M=\sum_{n\geq 1}\sum_{k=1}^{M-1}\left[
\frac{\mu\Big((Mn-k)^{\pm 2\pi i/\ln M}\Big)}{Mn-k}
f\Big(g\Big((Mn-k)^{\pm 2\pi i/\ln M}\Big)\Big)-\right.
\end{equation}
$$
\left.-\frac{\mu\Big((Mn)^{\pm 2\pi i/\ln M}\Big)}
{Mn}f\Big(g\Big((Mn)^{\pm 2\pi i/\ln M}\Big)\Big)\right].
 $$
\end{thrm}
\textbf{Доказательство.} Применим (\ref{eq-5}) к функции $F(z)=\mu(z)f(g(z)).\Box$

Большое семейство наборов точек для граничных представлений можно получить, применив к последовательности
$(Mn)^{2\pi i/\ln M}$ (или $(Mn)^{-2\pi i/\ln M}$) сюръективные функции $g\in\mathbf{A_B},$ важными примерами
которых служат произведения Бляшке \cite{M}. Все это также позволяет получать новые формулы суммирования для
интегралов по единичной окружности.

\begin{thrm}\label{Th-7}
Пусть целые $L\neq 0$  и $M>1,$ $\abs{a}=1,$ $b_1,\ldots,b_J$ по абсолютной величине меньше единицы, и
$f(z)\in\mathbf{A_B}.$ Если
$$
F(z)=f\left(\prod_{j=1}^J\frac{z+b_j}{1+\overline{b_j}\,z}\right)\quad
\mbox{и}\quad \alpha\prod_{j=1}^Jb_j=c,\quad\mbox{то}
 $$
$$
\ln M\,f(c)=\sum_{n\geq 1}\sum_{k=1}^{M-1}\left[\frac{1}{Mn-k}F\Big(a(Mn-k)^{2\pi iL/\ln M}\Big)
-\frac{1}{Mn}F\Big(a(Mn)^{2\pi iL/\ln M}\Big)\right].\Box
 $$
\end{thrm}

\begin{rmrq}
Другие весьма разнообразные выражения для сумм рядов и граничные представления значений функций рядами
связаны с представлениями (\ref{eq-22}) и (\ref{eq-23}), а также с формулами Пуассона, Шварца и Йенсена \cite{M}.
\end{rmrq}

\subsection{Логарифмы с рациональными основаниями}
\hspace{13pt} Здесь будут получены обобщения Теоремы \ref{Th-2}, которые позволят рассматривать в
качестве исходных последовательности узлов суммирования вида $\{\log_An\lmid n\in\mathds{N}\}$ с рациональными $A>1.$
\begin{thrm}\label{Th-4}
Если $M,N$ целые $>1,$ причем $M\neq N,$ и $f(x)\in\mathbf{P_B},$ то
\begin{equation}\label{eq-11}
\ln(M/N)\int_0^1f(t)dt=\sum_{n\geq 1}
\sum_{k=1}^{M-1}\left(\frac{1}{Mn-k}f\left(\frac{\ln(Mn-k)}{\ln(M/N)}\right)-\frac{1}{Mn}
f\left(\frac{\ln(Mn)}{\ln(M/N)}\right)\right)-
\end{equation}
$$
-\sum_{n\geq 1}\sum_{k=1}^{N-1}\left(\frac{1}{Nn-k}f\left(\frac{\ln(Nn-k)}{\ln(M/N)}\right)-\frac{1}{Nn}
f\left(\frac{\ln(Nn)}{\ln(M/N)}\right)\right).
 $$

Каждый из рядов по $n$ сходится абсолютно.
\end{thrm}
\textbf{Доказательство.} При $\Re(s)>0$ из (\ref{Hpt}) следует равенство
\begin{equation}\label{eq-2}
\left(N^{1-s}-M^{1-s}\right)\zeta(s)=\sum_{n\geq 1}\left[\sum_{k=1}^{M-1}\left(
\frac{1}{(Mn-k)^s}-\frac{1}{(Mn)^s}\right)-\sum_{k=1}^{N-1}\left(
\frac{1}{(Nn-k)^s}-\frac{1}{(Nn)^s}\right)\right].
\end{equation}

Поскольку $N^{1-s}-M^{1-s}=0,$ если $s=1+2\pi il/\ln(M/N),$ $l\in\mathds{Z},$ то при целых $l$ верно
$$
\sum_{n\geq 1}\Bigg[\sum_{k=1}^{M-1}\left(
\frac{1}{(Mn-k)^{1+2\pi il/\ln(M/N)}}-\frac{1}{(Mn)^{1+2\pi il/\ln(M/N)}}\right)-
 $$
$$
-\sum_{k=1}^{N-1}\left(
\frac{1}{(Nn-k)^{1+2\pi il/\ln(M/N)}}-\frac{1}{(Nn)^{1+2\pi il/\ln(M/N)}}\right)\Bigg]=\ln(M/N)\delta_{l,0}.
 $$
Дальнейшее доказательство проводится также как и в Теореме \ref{Th-2}. Сходимость и полиномиальная
аппроксимация исследуются по отдельности для рядов
$$
\sum_{n\geq 1}\left(\frac{1}{Mn-k}f\left(\frac{\ln(Mn-k)}{\ln(M/N)}\right)-\frac{1}{Mn}
f\left(\frac{\ln(Mn)}{\ln(M/N)}\right)\right),\;\; k=1\div M-1,\; \mbox{и}
 $$
$$
\sum_{n\geq 1}\left(\frac{1}{Nn-k}f\left(\frac{\ln(Nn-k)}{\ln(M/N)}\right)-\frac{1}{Nn}
f\left(\frac{\ln(Nn)}{\ln(M/N)}\right)\right),\;\; k=1\div N-1.\Box
 $$

\begin{rmrq}
Замена пары $(M,N)$ на ее целую кратную $(LM,LN),$ $L>1,$ приводит к иному ряду (\ref{eq-11}) на  более
разреженной подпоследовательности узлов, а ряд для $(LM,L)$ отличается от ряда в правой части (\ref{eq-4}), поскольку
с ростом $L$ появляется и растет разреженность ряда $\sum_{n\geq 1}C_n(L)f(\log_{(LM)/L}n).$ Т.е. для
каждого положительного рационального $A\neq 1$ с помощью линейных комбинаций мы можем получить бесчисленное
множество различных представлений интегралов по периоду в виде рядов по взвешенным значениям интегрируемой
функции в узлах $\{\Big\{\log_An\Big\}\lmid n\in\mathds{N}\}.$

Похожее разрежение ряда вызывается использованием параметра $L$ в (\ref{eq-5}).
\end{rmrq}

\begin{rmrq}
Подобно тому как Теоремы \ref{Th-25}-\ref{Th-3} раскрывают содержание Теоремы \ref{Th-2}, можно развить следствия Теоремы \ref{Th-4}. Вследствие чего полученные ранее свойства последовательностей
величин $\log _Mn$ при натуральных $M$ распространяются  на последовательности величин $\log _An$ с любым рациональным $A>1.$
\end{rmrq}

\subsection{Ряды связанные с производными $\zeta-$функции}\label{der} 
\hspace{13pt} Здесь предлагается еще один способ получения представлений интегралов по периоду рядами.

Пусть $L,M,N$ --- целые числа $>1$ и  $M\neq N.$ При $\Re{s>0}$ справедливы равенства (см. (\ref{eq-2}))
$$
\left(N^{1-s}-M^{1-s}\right)^L\zeta(s)=\sum_{l=0}^{L-1}C_{L-1}^l(-M^{1-s})^lN^{(L-1-l)(1-s)}\,\cdot\,
\left(N^{1-s}-M^{1-s}\right)\zeta(s)=
 $$
$$
=\sum_{l=0}^{L-1}C_{L-1}^l(-M)^lN^{L-1-l}\left[\sum_{k=1}^{M-1}\sum_{n\geq 1}
\Big((M^lN^{L-1-l}(Mn-k))^{-s}-(M^{l+1}N^{L-1-l}n)^{-s}\Big)-\right. 
 $$
\begin{equation}\label{eq-24}
\left. -\sum_{k=1}^{N-1}\sum_{n\geq 1}\Big((M^lN^{L-1-l}(Nn-k))^{-s}-(M^lN^{L-l}n)^{-s}\Big)\right].
\end{equation}

Продифференцируем первое и последнее выражение в (\ref{eq-24}) $L-1$ раз по $s.$ Положив
$s=1+2\pi ij/\ln{(M/N)},\; j\in\mathds{Z},$ получаем равенство
$$
(L-1)!\,\Big(\ln(M/N)\Big)^L\,\delta_{j,0}=\sum_{l=0}^{L-1}C_{L-1}^l(-M)^lN^{L-1-l}\;\times
 $$
$$
\times\left[\sum_{k=1}^{M-1}\sum_{n\geq 1}\Bigg(\frac{(-\ln(M^lN^{L-1-l}(Mn-k)))^{L-1}}{(M^lN^{L-1-l}(Mn-k))^s}-
\frac{(-\ln(M^{l+1}N^{L-1-l}))^{L-1}}{(M^{l+1}N^{L-1-l}n)^s}\Bigg)-\right. 
 $$
$$
\left. -\sum_{k=1}^{N-1}\sum_{n\geq 1}\Bigg(\frac{(-\ln(M^lN^{L-1-l}(Nn-k)))^{L-1}}{(M^lN^{L-1-l}(Nn-k))^s}-
\frac{(-\ln(M^lN^{L-l}))^{L-1}}{(M^lN^{L-l}n)^s}\Bigg)\right]\Bigg|_{s=1+2\pi ij/\ln{(M/N)}},
 $$
что позволяет сформулировать следующее утверждение.
\begin{thrm}\label{Th-44}
Если $L,M,N$ целые $>1,$ $M\neq N,$ и $f(x)\in\mathbf{P_B}_{L-1},$ то
\begin{equation}\label{eq-44}
(L-1)!\,\Big(\ln(M/N)\Big)^L\int_0^1f(t)dt\;=\;\sum_{l=0}^{L-1}C_{L-1}^l(-M)^lN^{L-1-l}\;\times
\end{equation}
$$
\times\left[\sum_{k=1}^{M-1}\sum_{n\geq 1}\Bigg(\,f\Big(M^lN^{L-1-l}(Mn-k),L-1\Big)-
\,f\Big(M^{l+1}N^{L-1-l}n,L-1\Big)\Bigg)-\right. 
 $$
$$
\left. -\sum_{k=1}^{N-1}\sum_{n\geq 1}\Bigg(\,f\Big(M^lN^{L-1-l}(Nn-k),L-1\Big)-
\,f(M^lN^{L-l}n,L-1)\Bigg)\right],
 $$
где $f(w,L-1)=(-\ln w)^{L-1}w^{-1}f(\log_{M/N}w).$

Каждый из рядов по $(n,k)$ сходится абсолютно.
\end{thrm}
\textbf{Доказательство} необходимой нам сходимости проводим по использованной выше схеме, т.е. отдельно для
каждого ряда по $n\geq 1$ при фиксированной паре $(l,k).$ Множители ${\ln(w)}^{L-1}$ указывают на необходимость
использования других классов функций$.\Box$

Теоремы \ref{Th-25}-\ref{Th-3}  очевидным образом переносятся и на этот случай.

\begin{rmrq}
Теорема \ref{Th-44} останется справедливой и при $N=1,$ $M>1,$ если из правой части равенства (\ref{eq-44}) убрать
сумму, начинающуюся с $\sum_{k=1}^{N-1}\sum_{n\geq 1}.$
\end{rmrq}

\section{Применения интегральных представлений $\zeta-$функции}
\hspace{13pt} Ниже будут предложены применения тех же формул Эйлера и некоторых из использованных выше идей к
интегральным представлениям $\zeta-$функции.

Для целых $J\geq 0,$ $M>1,$ положительного $a$ и $b,$ $\Re(b)>-1,$ определим функцию
\begin{equation}\label{eq-9}
\Phi_{a,b,M,J}(w)=(-1)^J\frac{d^J}{dw^J}\,\sum_{n\geq 1}
\sum_{k=1}^{M-1}\left[(Mn-k)^be^{-w(Mn-k)^a}-(Mn)^be^{-w(Mn)^a}\right].
\end{equation}

\begin{thrm}\label{Th-2-int}
Пусть $z\in\mathds{R},$ целые $M\geq 2,$ $J\geq 0,$ $a>0,$ и $\Re(b)>-Ja-1.$ Справедливо равенство
\begin{equation}\label{eq-10}
\sum_{n\in\mathds{Z}}M^{(Ja+1+b)(n+z)}\Phi_{a,b,M,J}\left(M^{a(n+z)}\right)=
\Gamma\Big((1+b)/a+J\Big)/a.
\end{equation}
\end{thrm}
\textbf{Доказательство.} Пользуясь соотношением
$n^{-s}=n^b/\Gamma\Big((s+b)/a\Big)\cdot\int_0^\infty x^{(s+b)/a-1}e^{-xn^a}dx,$ при $\Re{(s+b)}>0$ перепишем
правую  часть (\ref{Hpt}) в интегральной форме
$$
\Gamma\Big((s+b)/a\Big)^{-1}\int_0^{\infty}dt\cdot t^{(s+b)/a-1}
\sum_{n\geq 1}\sum_{k=1}^{M-1}\left((Mn-k)^be^{-(Mn-k)^at}-(Mn)^be^{-(Mn)^at}\right),
 $$
из которой $J$ интегрированиями по частям
$\Big[\int_0^\infty t^{\alpha}f(t)dt=-(\alpha+1)\int_0^\infty t^{\alpha+1}f'(t)dt\Big]$ выводим:
$$
\left(1-M^{1-s}\right)\zeta(s)=\frac{1}{\Gamma\Big((s+b)/a+J\Big)}\int_0^{\infty}dt\cdot t^{(s+b)/a+J-1}
\Phi_{a,b,M,J}(t).
 $$

Это соотношение аналитически продолжается на значения $s,b,$ где $\Re{(s+b)/a+J>0}.$

Подставляя значения $s=1-2\pi il/\ln M,$ $l\in\mathds{Z},$ получаем равенства
$$
\Gamma\Big((1+b)/a+J\Big)\,\delta_{l,0}\ln M=\int_0^{\infty}dt\cdot t^{J+(1+b)/a-1-2\pi il/(a\ln M)}
\Phi_{a,b,M,J}(t)=
 $$
$$
=a\ln M\int_{-\infty}^{\infty}dw\cdot e^{w[(Ja+1+b)\ln M-2\pi il]}\Phi_{a,b,M,J}(e^{w\,a\ln M})=
 $$
$$
=a\ln M\int_0^1dw\cdot e^{-2\pi ilw}\sum_{n\in\mathds{Z}}M^{(w+n)[(Ja+1+b)]}\Phi_{a,b,M,J}
(M^{(w+n)a}),
 $$
из которых следует, что функция в последнем интеграле не зависит от $w.\Box$

Наш следующий результат использует еще одно соотношение для $\zeta-$функции.

\begin{thrm}\label{Th-4-int}
Пусть $M,N\in\mathds{N},$ причем $M>1<N$ и $M\neq N,$ целое $J\geq 0,$ $a>0,$ $\Re(b)>-Ja-1,$ $z\in\mathds{R}$ и
функция $\Phi_{a,b,M,J}(w)$ определена соотношением (\ref{eq-9}). Справедливо равенство
\begin{equation}\label{eq-8}
\sum_{n\in\mathds{Z}}\left(\frac{M}{N}\right)^{(aJ+1+b)(n+z)}
\left[\Phi_{a,b,M,J}\left(\left(\frac{M}{N}\right)^{a(z+n)}\right)-
\Phi_{a,b,N,J}\left(\left(\frac{M}{N}\right)^{a(z+n)}\right)\right]=\frac{\Gamma(\frac{1+b}{a}+J)}{a}.
\end{equation}
\end{thrm}
\textbf{Доказательство.} Следуем общему плану доказательства Теоремы \ref{Th-2-int}.

Преобразуем (\ref{eq-2}) к виду
$$
\left(N^{1-s}-M^{1-s}\right)\zeta(s)=
\Gamma((s+b)/a)^{-1}\int_0^{\infty}dt\cdot t^{(s+b)/a-1}\cdot
 $$
$$
\cdot\sum_{n\geq 1}\left(\sum_{k=1}^{M-1}\Big[(Mn-k)^be^{-t(Mn-k)^a}-(Mn)^be^{-t(Mn)^a}\Big]-\right.
 $$
$$
\left.-\sum_{k=1}^{N-1}\Big[(Nn-k)^be^{-t(Nn-k)^a}-(Nn)^be^{-t(Nn)^a}\Big]\right)=
 $$
(после интегрирований по частям)
$$
=\frac{1}{\Gamma\Big((s+b)/a+J\Big)}
\int_0^{\infty}dt\cdot t^{(s+b)/a-1+J}
\left(\Phi_{a,b,M,J}(t)-\Phi_{a,b,N,J}(t)\right).
 $$
Подставляя значения $s=1-2\pi il/\ln(M/N),$ $l\in\mathds{Z},$ получаем равенства
$$
\Gamma\Big(\frac{1+b}{a}+J\Big)\,\delta_{l,0}=\frac{1}{\ln(M/N)}\int_0^{\infty}dt\cdot
t^{J+(1+b)/a-1-2\pi il/(a\ln(M/N))}\Big(\Phi_{a,b,M,J}(t)-\Phi_{a,b,N,J}(t)\Big)=
 $$
$$
=a\int_{-\infty}^{\infty}dw\cdot e^{w[(Ja+1+b)\ln(M/N)-2\pi il]}
\Big(\Phi_{a,b,M,J}(e^{wa\ln(M/N)})-\Phi_{a,b,N,J}(e^{wa\ln(M/N)})\Big)=
 $$
$$
a\int_0^1dw\cdot e^{-2\pi ilw}\sum_{n\in\mathds{Z}}\left(\frac{M}{N}\right)^{(aJ+1+b)(n+z)}
\left[\Phi_{a,b,M,J}\left(\left(\frac{M}{N}\right)^{a(z+n)}\right)-
\Phi_{a,b,N,J}\left(\left(\frac{M}{N}\right)^{a(z+n)}\right)\right],
 $$
откуда следует наше утверждение$.\Box$

\begin{bspl}
Если целое $b\geq 0,$ то $\Phi_{1,b,M,J}(w)=(-1)^{J+b}\frac{d^{J+b}}{dw^{J+b}}\,
\left(\frac{1}{e^w-1}-\frac{M}{e^{wM}-1}\right),$
из чего следуют равенства
$$
\sum_{n\in\mathds{Z}}\frac{2^{n+z}}{e^{2^{n+z}}+1}=1,\;\;
\sum_{n\in\mathds{Z}}\frac{4^{n+z}e^{2^{n+z}}}{\left(e^{2^{n+z}}+1\right)^2}=1,\;\;
\sum_{n\in\mathds{Z}}\frac{3^{n+z}\left(e^{3^{n+z}}+2\right)}{e^{2\cdot 3^{n+z}}+e^{3^{n+z}}+1}=1,
 $$
$$
\sum_{n\in\mathds{Z}}\left(\frac{3}{2}\right)^{n+z}\frac{e^{\left(\frac{3}{2}\right)^{n+z}}+1}
{e^{3\left(\frac{3}{2}\right)^{n+z}}+2e^{2\left(\frac{3}{2}\right)^{n+z}}+2e^{\left(\frac{3}{2}\right)^{n+z}}+1}=1.
 $$

Очевидно, что каждое из приведенных соотношений справедливо при всех $z,$ для которых не обращаются в ноль
знаменатели членов ряда в этом соотношении.
\end{bspl}

\begin{bspl}\label{bspl} 
При $\Re(b)>0$ функцию $\Phi_{a,b,M,0}$ можно переписать в виде
$$
\Phi_{a,b,M,0}(w)=\Big(\theta_{a,b}(w)-M\,\theta_{a,b}(wM^a)\Big)/2,\quad
\theta_{a,b}(w)=\sum_n\abs{n}^be^{-w\,\abs{n}^a}.
 $$
В частности $\theta_{2,0}(w)$ с точностью до обозначений аргументов --- $\theta-$функция Якоби.
\end{bspl}

\begin{bspl}
При $b=\beta a,$ где целое $\beta\geq 0,$ с использованием $\theta$ из Примера \ref{bspl} получаем
$$
\Phi_{a,\beta a,M,J}(w)=\frac{(-1)^{J+\beta}}{2}\frac{d^{J+\beta}}{dw^{J+\beta}}
\Big(\theta_{a,0}(w)-M\theta_{a,0}(wM^a)-1+M\Big).
 $$
\end{bspl}

Отправляясь от представлений производных $\zeta-$функции, можно получить еще одно 
семейство рядов.
\begin{thrm}\label{Th-5-int}
Пусть целые $J,M,N,L,j_1,l_1,\ldots$ неотрицательны, причем $M\neq N,$ $N>1<L$ и $j_1+j_2+\ldots=J,$
$l_1+l_2+\ldots=R-1,$ а также $a>0,$ $\Re(b)>-Ja-1,$ $z\in\mathds{R}.$

В случае функции
$$
\Psi_{a,b,M,N,J,R,\textbf{j},\textbf{r}}(t)=(-1)^J\frac{d^{j_1}}{dt^{j_1}}\Bigg(\Bigg(\frac{\ln t}{a}\Bigg)^{l_1}\cdot
\frac{d^{j_2}}{dt^{j_2}}\Bigg(\Bigg(\frac{\ln t}{a}\Bigg)^{l_2}\cdot\ldots\sum_{l=0}^{L-1}C_{L-1}^l(-M)^lN^{L-1-l}
\,\times
 $$
$$
\times\left[\sum_{k=1}^{M-1}\sum_{n\geq 1}\!
\Big((M^lN^{L-1-l}(Mn-k))^be^{-t(M^lN^{L-1-l}(Mn-k))^a}\!-\!(M^{l+1}N^{L-1-l}n)^be^{-t(M^{l+1}N^{L-1-l})^a}\Big)
-\right. 
 $$
$$
\left. -\sum_{k=1}^{N-1}\sum_{n\geq 1}\Big((M^lN^{L-1-l}(Nn-k))^be^{-t(M^lN^{L-1-l}(Nn-k))^a}-(M^lN^{L-l}n)^b
e^{-t(M^lN^{L-l})^a}\Big)\right]
 $$

справедливо равенство
$$
\sum_{n\in\mathds{Z}}\Psi_{a,b,M,N,J,R,\textbf{j},\textbf{l}}\left(\left(\frac{M}{N}\right)^{a(z+n)}\right)=
(L-1)!\;\Gamma((1+b)/a+J)/a.
 $$
\end{thrm}
\textbf{Доказательство.} Равенство (\ref{eq-24}) можно преобразовать к виду
$$
\left(N^{1-s}-M^{1-s}\right)^L\zeta(s)=\Gamma\Big(\frac{s+b}{a}\Big)^{-1}\int_0^{\infty}dt
\cdot t^{(s+b)/a-1}\cdot \sum_{l=0}^{L-1}C_{L-1}^l(-M)^lN^{L-1-l}\,\times
 $$
$$
\times\left[\sum_{k=1}^{M-1}\sum_{n\geq 1}\!
\Big((M^lN^{L-1-l}(Mn-k))^be^{-t(M^lN^{L-1-l}(Mn-k))^a}\!-\!(M^{l+1}N^{L-1-l}n)^be^{-t(M^{l+1}N^{L-1-l})^a}\Big)
-\right. 
 $$
$$
\left. -\sum_{k=1}^{N-1}\sum_{n\geq 1}\Big((M^lN^{L-1-l}(Nn-k))^be^{-t(M^lN^{L-1-l}(Nn-k))^a}-(M^lN^{L-l}n)^b
e^{-t(M^lN^{L-l})^a}\Big)\right].
 $$
 
Продифференцируем $L-1$ раз обе части последнего равенства. При этом в правой части под знаком хорошо сходящегося
интеграла чередуем дифференцирования по $s$ с интегрированиями по частям вида
$\int_0^\infty t^{\alpha}f(t)dt=-(\alpha+1)\int_0^\infty t^{\alpha+1}f'(t)dt.$
$$
d^{L-1}\Big[\left(N^{1-s}-M^{1-s}\right)^L\zeta(s)\Big]\Big/ds^{L-1}=
 $$
$$
=\frac{\left(\ln M/N\right)^{L-1}}{\Gamma\Big((s+b)/a+J\Big)}\int_0^{\infty}dt\cdot t^{(s+b)/a-1+J}
\Psi_{a,b,M,N,J,R,\textbf{j},\textbf{r}}(t).
 $$
Далее также как и в двух предшествующих доказательствах подставим вместо $s$ величины
$s=1-2\pi il/\ln(M/N),$ $l\in\mathds{Z},$ выполним замену переменных $t=e^{w\;2\pi ia\ln M/N}$ и разделим область
интегрирования на отрезки с концами в целых числах$.\Box$

\begin{rmrq}
Как и в случае Теоремы \ref{Th-44} доказанная Теорема \ref{Th-5-int} останется справедливой и при $N=1,$ $M>1,$ если
из определения функции $\Psi_{a,b,M,N,J,R,\textbf{j},\textbf{r}}(t)$ убрать сумму, начинающуюся с
$\sum_{k=1}^{N-1}\sum_{n\geq 1}.$
\end{rmrq}

\begin{rmrq}
В доказательствах этого раздела на параметры $a,b,z\ldots$ были наложены ограничения, которые могут ослаблены в
силу аналитичности полученных соотношений.
\end{rmrq}

\section{Обобщения}

\subsection{}
\hspace{13pt} В доказательствах Теорем \ref{Th-2} и \ref{Th-4} мы рассматривали представления интегралов в виде
конечных сумм (по $1\leq k<M,N$) более простых рядов и приводили доказательства при условиях на функции, не
зависящих от $M$  или $N.$ Не разбивая ряды, мы могли бы расширить классы функций $\mathbf{P_B}$  и
$\mathbf{C_B},$ для которых
справедливы полученные разложения.

\subsection{}
\hspace{13pt} Здесь рассматриваются различные функции, которые способны в большей или меньшей степени заменить
$\zeta-$функцию при получении сходных с приведенными ранее рядами.

Большое число разложений в ряды связано с функциями вида
\begin{equation}\label{eq-15}
\widehat{\zeta}(s)=\sum_{n_1,\ldots,n_d\geq 1}\Psi(n_1,\ldots,n_d;s),
\end{equation}
если при натуральных $M$ справедливо
$\Psi(Mn_1,\ldots,Mn_d;s)=M^{-s}\Psi(n_1,\ldots,n_d;s).$

В качестве примеров укажем две аналитические функции, которые при $\Re{(s)}>d$ определяются как
$$
\widehat{\zeta}_1(s)=\sum_{n_1>1,\ldots,n_d\geq 1}\frac{\Big(\sum_{1\leq j\leq d}{n_j}\Big)^{s+\pi}}
{\sum_{k=1}^K\Big(\sum_{1\leq j\leq d}A_{kj}{n_j}\Big)^{2s+\pi}+\left(\prod_{ j=1}^{d}n_j\right)^{(2s+\pi)/d}}\,,
 $$
$$
\widehat{\zeta}_2(s)=\sum_{n_1>1,\ldots,n_d\geq 1}\frac{\sum_{1\leq  j\leq d-1}jn_jn_{j+1}}{\left(\sum_{1\leq  j
\leq d}{n_j}\right)^2}\,\cdot\,\frac{1}{\prod_{i=1}^I\left(\sum_{ j=1}^{d}B_{ij}n_j\right)^{\upsilon_is}}\,,
 $$
где все величины $A_{kj},B_{ij}$ положительны, и $\sum_i\upsilon_i=1.$

\begin{thrm}
Функции $\widehat{\zeta}_{1,2}(s)$ обладают свойствами: 

1) если $\Re{(s)}>d-1$ и целое $M>1,$ то
\begin{equation}
(1-M^{d-s})\widehat{\zeta}(s)=
\end{equation}
$$
=\sum_{n_1,\ldots,n_d\geq 1}\Bigg(\sum_{0\leq k_1,\ldots,k_d\leq M-1}\Psi(Mn_1-k_1,\ldots,Mn_d-k_d;s)-M^d
\Psi(Mn_1,\ldots,Mn_d;s)\Bigg),
 $$
причем при всяком $\varepsilon>0$ сходимость ряда по $(n_1,\ldots,n_d)$ равномерна в области
$\Re{(s)}>d-1+\varepsilon;$

2) у функций в области $\Re{(s)}>d-1$ лишь один полюс ($s=d$);

3 ) при всех целых $M>1$ величина
$$
(\ln M)^{-1}\sum_{n_1,\ldots,n_d\geq 1}\Bigg(\sum_{0\leq k_1,\ldots,k_d\leq M-1}
\Psi(Mn_1-k_1,\ldots,Mn_d-k_d;d)-M^d\Psi(Mn_1,\ldots,Mn_d;d)\Bigg)
 $$
для каждой из наших функций одинакова и положительна$.\Box$
\end{thrm}

Эти свойства позволяют получить полные аналоги Теорем (\ref{Th-2})-(\ref{Th-44}). 

Для получения аналогов Теорем (\ref{Th-2-int}), (\ref{Th-4-int}) и  (\ref{Th-5-int}) подходит функция
$\widehat{\zeta}_2$ равно как и многие другие функции вида
$$
\widehat{\zeta}_3=\sum_{n_1,\ldots,n_l\geq 1}\xi(n_1,\ldots,n_l)\,\psi(n_1,\ldots,n_l)^{-s},
 $$
в которых $\xi$ однородная функция нулевого, а $\psi$ --- первого порядка. Это связано с равенством
$$
(1-M^{d-s})\widehat{\zeta}_3(s)\,=\,\frac{(-1)^J}{\Gamma(s/\alpha+J)}\int_0^{\infty}dt\cdot
t^{s/\alpha+J-1}\cdot
 $$
$$
\cdot\frac{d^J}{dt^J}\left(\sum_{0\leq k_1,\ldots,k_d\leq M-1}\xi(M\textbf{n}-\textbf{k})\,e^{-t\,\psi(M\textbf{n}-
\textbf{k})^\alpha}-M^d\xi(M\textbf{n})\,e^{-t\,\psi(M\textbf{n})^\alpha}\right).
 $$

Отметим возможность других подобных обобщений дзета-функции, связанных с использованием в соотношениях
подобных (\ref{eq-15}) других областей суммирования по $(n_1,\ldots,n_d).$ При этом полезно выполнение требований:
а)~область не содержит точки $(0,\ldots,0)$ и б)~область замкнута относительно сдвигов на векторы $(n_1,\ldots,n_d)$ с
целыми неотрицательными компонентами. Разработанная сегодня теория автоморфных форм может быть очено полезна
в этих рассмотрениях.

\subsection{}
\hspace{13pt} Сходные идеи могут использованы в связи с другими функциями.

Пусть  $a_1+a_2\neq b_1+b_2,$ $C(z)$ --- целая функция и $C(a)e^{ax}=\sum_{n\geq 0}C_n(x)a^n.$ Это
позволяет получить равенства
$$
C(x)\frac{e^{(a_1+a_2)x}-e^{(b_1+b_2)x}}{x}=
\sum_{n\geq 0}\frac{e^{a_2x}C_n(a_1)-e^{b_2x}C_n(b_1)}{x}x^n.
 $$
Подставим $x_l=2\pi il/(a_1+a_2-b_1-b_2)$
$$
C(0)(a_1-b_1)\delta_{l,0}=\sum_{n\geq 0}\left(e^{\frac{2\pi ila_2}{a_1+a_2-b_1-b_2}}C_{n+1}(a_1)-
e^{\frac{2\pi ilb_2}{a_1+a_2-b_1-b_2}}C_{n+1}(b_1)\right)\left(\frac{2\pi il}{a_1+a_2-b_1-b_2}\right)^{n+1}
 $$
и при $f(z)=\sum\limits_{l=L_1}^{L_2}f_le^{2\pi ilz}$ получим равенство
$$
C(0)(a_1-b_1)\int_0^1f(z)dz=\sum_{n\geq 0}\frac{1}{(a_1+a_2-b_1-b_2)^{n+1}}\times
 $$
$$
\times\frac{d^{n+1}}{dz^{n+1}}
\Bigg(f\left(z+\frac{a_2}{a_1+a_2-b_1-b_2}\right) C_{n+1}(a_1)-
f\left(z+\frac{b_2}{a_1+a_2-b_1-b_2}\right)C_{n+1}(b_1)\Bigg)\Bigg|_{z=0}.
 $$

Последнее равенство можно распространить на функции $f(z)$ с периодом единица и голоморфные в
подходящей полосе.

\section{}
\hspace{13pt} Автор признателен Ф.С.Джепарову, Д.В.Львову и А.Н.Тюлюсову за очень полезное обсуждение  работы.

Victor.Shestopal@itep.ru, veshestopal1@yandex.ru

\clearpage
\end{document}